\documentstyle[12pt]{article}
\newcommand{\R}{{ R\hspace*{-0.9ex}\rule{0.15ex}{1.5ex}\hspace*{0.9ex}}}
\newcommand{\C}{I\!\!\!\!C}
\def\N{I\!\!N}
\newtheorem{theo}{Th\'eor\`eme}[section]
\newtheorem{deff}{D\'efinition}[section]
\newtheorem{prop}{Proposition}[section]
\newtheorem{lem}{Lemme}[section]

\newtheorem{coro}{Corollaire}[section]

\begin{document}
\title{Forme Equivalente \`a la Condition $\Delta_2 $ et Certains r\'esultats de s\'eparations dans les Espaces Modulaires  }
\date{}
\author{Par \\ A.Hajji }
\maketitle
{\bf Abstract}. In this paper, we present an equivalent form of the  $\Delta_2 $-condition which allow us to redefine the topological vector space structure of a modular spaces using the \\ filter base. We show also the characterization of closed subsets (in the sens of this topology ) of a modular spaces which permit us to establish some separation results in modular spaces.\\

A.M.S Subject Classifications: 46A80.
\section{Introduction}

\quad La condition $\Delta_2$  joue un r\^ole primordial dans la th\'eorie des espaces modulaires. J.Musielak \cite{jm} a montr\'e qu' une base modulaire est topologique si et seulement si le modulaire $\rho$ v\'erifie la condition $\Delta_2$, et en appliquant, pour un espace modulaire $X_\rho $ avec le modulaire $\rho$, deux formes \'equivalentes \`a la condition $\Delta_2$\\
 (1) Pour toute suite $(x_n )_{n\in \N}$ dans $X_\rho $ si $\lim\limits_{n\to +\infty}\rho(x_n)=0$ alors
$\lim\limits_{n\to +\infty}\rho(2x_n)=0.$\\
(2) Pour tout $u\in {\cal{B}}$ il existe $v\in {\cal{B}}$ tel que $2v\in u$, \\
o\`u  ${\cal{B}} $ est la famille de $\rho$-boules $B_\rho (0,r) =\{ x\in X\rho: \rho (x) <r \}$ pour $r>0$. La forme (1) de la condition $\Delta_2 $ est \'etroitement li\'ee avec l'espace modulaire $X_\rho $, la forme (2) de la condition $\Delta_2$ se r\'ef\`ere \`a la topologie (\`a base ${\cal{B}}$ ) introduite dans $X_\rho $.\\

En g\'en\'eral, la $\rho$-convergence (convergence au sens de $\rho$: pour toute suite $(x_n)_{n\in \N}$  dans $X_\rho , \ \  x_n \stackrel{\rho}{\rightarrow} x $ est \'equivalent qu'il existe $\lambda >0$ tel que la suite $(\rho(\lambda (x_n -x))_{n\in \N}$ converge vers z\'ero)  n 'implique pas  la convergence en norme, tendis que si $\rho $ est convexe et v\'erifie la condition $\Delta_2 $, alors  ces  deux modes de convergences  sont \'equivalents ( voir  \cite{m1},  \cite{krk}) et \cite{jm}).\\

 Nous pr\'esentons dans ce papier,  la suivante forme \'equivalente \`a la condition $\Delta_2$  (voir le Lemme 1.1 ) \\
(3) Pour tout $\epsilon >0$, il existe $L>0,$ et il existe $\delta>0$ tels que si $\rho (x) <L $ et $\rho (x-y)< \delta$ alors $|\rho (y) -\rho(x)|<\epsilon $.\\
C'est \`a dire que le modulaire $\rho$ poss\`ede la propri\'et\'e d'une " certaine" continuit\'e uniforme local. La forme (3) de la condition $\Delta_2 $, nous permit de red\'efinir, \`a l' aide  d'une base de filtres (voir \cite{jd} ), la structure d'espace vectoriel topologique de l'espace modulaire $X_\rho $ et de munir l'espace moduliare $X_\rho $ d'une topologie s\'epar\'ee not\'ee $\tau$, cette derni\`ere propri\'et\'e pour la topologie $\tau$ n'a pas \'et\'e mentionn\'ee dans \cite{jm}. Aussi la forme (3) de la condition $\Delta_2$ est fr\'equement utiliser pour caracteriser les ensembles ferm\'es au sens de la topologie $\tau$, ce qui donne la d\'emonstration de la d\'efinition 5.4 \cite{jm} des ensembles $\rho$-ferm\'es dans l'espaces modulaire $X_\rho,$ et pour montrer certains r\'esultats de s\'eparations dans l'espace modulaire $X_\rho$.  \\

Nous commen\c{c}ons par un  rappel sur les notions de base de la th\'eorie des espaces modulaires ( voir \cite{jm}).
  \begin{deff} Soit $X$ un espace vectoriel sur $\R$ ou $\C$.

\quad a) La fonctionnelle $\rho: X\longrightarrow [0,+\infty]$ est
dite un modulaire si pour tout $x,y \in X$ on a:

\qquad i) $\rho(x)=0$ si et seulement si $ x = 0$.

\qquad ii) $\rho(-x)=\rho(x)$ dans le cas r\'eel

\qquad et  $\rho(e^{it}x)=\rho(x)$ pour tout $t\in \R$ dans le cas
complexe.

\qquad iii) $\rho(\alpha x+\beta y)\leq \rho(x)+\rho(y)$ pour
$\alpha,\beta \geq 0$, $\alpha+\beta =1$, $\forall x, y\in X$.

\quad b) L'espace modulaire $X_\rho$ associ\'e au modulaire
$\rho$ est donn\'e par: $$ X_\rho=\{ x\in X \,  / \, \rho(\lambda
x)\rightarrow 0 \,\mbox{ \, quand \,} \lambda \rightarrow0 \}.$$
\end{deff}
{\bf Remarques }\\

\quad 1) Si on remplace iii) par iii$^{'}$):
  $$\rho(\alpha x+\beta
y)\leq \alpha^s\rho(x)+\beta^s\rho(y)$$
 pour $\alpha,\beta \geq 0$, $\alpha^s +\beta^s =1$
avec un $s\in ]0,1]$, alors le modulaire $\rho$ est dit un
s-convexe. \\
Un modulaire 1-convexe est dit convexe. \\

\quad 2) Si $\rho$ est un modulaire convexe alors l'espace modulaire $X_\rho $ sera \'equip\'e de la norme de Luxemburg d\'efinie par:
$$ \|x \|_{\rho}  = Inf \{u >0 , \rho( \frac{x}{ u  }) \leq 1 \} $$
Ou par la norme ( voir \cite{ab}), o\`u $\rho$ est suppos\'e convexe, d\'efinie par:
$$ \|f\|= \frac{\xi}{s(f)},   $$
o\`u $\xi >0$ un nombre reel, $f\in L^\psi$ et $s(f) = sup \{ s: \rho (sf) \leq \xi  \}>0.$\\
Par cons\'equent, si $\rho$ n'est pas convexe alors on ne peut pas conclure que $X_\rho$ est un espace norm\'e.\\

3) Si $\rho$ v\'erifie la forme (1) de la condition $\Delta_2 $, alors on a l'\'equivalence suivante:
$$  x_n  \to x \Longleftrightarrow x_n \stackrel{\rho}{\rightarrow} x, $$
o\`u $x_n \to x $ veut dire $\lim\limits_{n\to +\infty}\rho(x_n -x)=0$. En effet:\\
 Soit $x_n \stackrel{\rho}{\rightarrow} x$ alors il existe $\lambda >0 $ tel que la suite $ (\rho (\lambda (x_n -x))_{n\in \N} $ converge vers z\'ero. \\
Si $\lambda \geq 1 $ alors on aura: $\rho (x_n -x ) \leq \rho (\lambda (x_n -x))$ et ceci montre que $\lim\limits_{n\to +\infty}\rho(x_n -x)=0$.\\
Si $\lambda <1 $ alors il existe tout jour $p\in \N $ ( $\N$: l'ensemble  des entiers naturels) tel que $1\leq 2^p \lambda $ et on aura $\rho (x_n -x ) \leq \rho ( 2^p \lambda (x_n -x))$, or par la forme (1) de la condition $\Delta_2 $ on a $\lim\limits_{n\to +\infty}\rho(2^p \lambda(x_n -x))=0$, ce qui montre que $x_n \stackrel{\rho}{\rightarrow} x$ implique $\lim\limits_{n\to +\infty}\rho(x_n -x)=0$.\\
Inversement.\\
Soit $ x_n \to x$  veut dire il existe $\lambda =1$ tel que $\lim\limits_{n\to +\infty}\rho(x_n -x)=0$. D'o\`u  $x_n \stackrel{\rho}{\rightarrow} x$, finalement on peut conclure que si $\rho$ v\'erifie la forme (1) de la condition $\Delta_2 $, alors  $ x_n  \to x \Longleftrightarrow x_n \stackrel{\rho}{\rightarrow} x$ \\

\quad Rappelons que nous ne pouvons appliquer la forme (1) de la condition $\Delta_2 $  que pour les suites. Dans le r\'esultat suivant nous montrons une forme \'equivalente de la condition $\Delta_2 $ qui est  une sorte de continuit\'e uniforme local de $\rho $ et qui sera fr\'equemment utilis\'ee par la suite. \\
On dit que $\rho$ v\'erifie la propri\'et\'e de Fatou si  $ x_n \stackrel{\rho}{\rightarrow} x $ et $ y_n \stackrel{\rho}{\rightarrow} y $, alors $ \rho (x-y) \leq \liminf \rho(x_n - y_n ) $.
\begin{lem}  Soit $\rho$ un modulaire de $X_\rho $ on a:\\
 $\rho$ v\'erifie la condition $\Delta_2$
si et seulement si pour tout $\epsilon >0$, ils existent $L>0$ et
$\delta >0$ tels que si $\rho(x)<L$ et $\rho(x-y)<\delta$, alors
$|\rho(y)-\rho(x)|<\epsilon$.
\end{lem}
 {\bf Preuve du Lemme 1.1.}

$\Leftarrow )$ Soit la suite $(x_n)_{n\in \N}$ dans $X_\rho$ tel
que $\rho(x_n)\longrightarrow 0$ quand $n\to +\infty$. Ceci impliquera que
 pour tout $\epsilon >0$, il existe $n_0$ tel que pour $n>n_0$ on a: $$\rho(x_n)<\inf (L,\delta, \frac \epsilon 2). $$
 On pose
$X_n =x_n$  et $Y_n =2x_n$, pour $n>n_0$ on a:\\
$$\rho(X_n)=\rho(x_n) =
\rho(Y_n-X_n)<\inf (L,\delta, \frac \epsilon 2),$$
 donc
$\rho(Y_n)=\rho(2x_n)\leq \frac \epsilon 2 +\rho(x_n)\leq
\epsilon$ d\`es que $n>n_0$.\\ D'o\`u la suite $(\rho(2x_n))_{n\in \N}$ tend vers
z\'ero quand $n$ tend vers $+\infty$. Par suite $\rho$ v\'erifie la
condition $\Delta_2$.\\

$\Rightarrow )$ Soit $\rho$ un modulaire qui v\'erifie la condition
$\Delta_2$.\\
 Supposons qu'il existe $ \alpha >0 $ tels que pour tout $ L>0 $ et tout $\delta >0$, il existe $ x,y $ v\'erifiant  $\rho(x)<L, \rho(x-y)<\delta$ et
$|\rho(y)-\rho(x)|\geq \alpha$.\\
 Posons $L=\delta=\frac 1n$ alors
il existe $ x_n,y_n$ tels que
$$\rho(x_n)<\frac 1n ,\  \rho(y_n-x_n)<\frac 1n \  \mbox{et} \  |\rho(y_n)-\rho(x_n)|\geq \alpha, $$
 ce qui entraine que
$\rho(x_n)\longrightarrow 0$ et $\rho(y_n-x_n) \longrightarrow 0$
quand $n\to +\infty$.\\
 Or
$$\rho(y_n)= \rho((x_n - y_n )+ y_n )\leq \rho(2(x_n-y_n))+
\rho(2x_n) ,$$
 et comme $\rho$ v\'erifie la condition $\Delta_2$,
alors $\rho(y_n) \longrightarrow 0$ quand $n\to +\infty$, donc $$ |\rho(y_n)-\rho(x_n)| \longrightarrow 0 \ \mbox{ quand}\ n\to +\infty ,$$
 absurde car  $|\rho(y_n)-\rho(x_n)|\geq \alpha >0$ pour tout $n\in \N$. D'o\`u
$$\forall  \epsilon >0, \exists L>0, \exists \delta
>0, \ \
 \mbox{tels que si} \ \
 \rho(x)<\delta \ \ \mbox{et} \ \ \rho(y-x)<\delta ,
 \ \ \mbox {alors}\ \
 |\rho(y)-\rho(x)|<\epsilon.$$
Ce qui termine la preuve du Lemme 1.1.
\section{ Structure d'espace vectoriel topologique s\'epar\'e  d'un espace modulaire.}

\quad On consid\`ere l'espace modulaire $X_\rho$ o\`u $\rho$ est un modulaire v\'erifiant la condition $\Delta_2$. Soit la famille
${\cal B}=\{B_\rho(0,r)\, /\, r>0\}$ o\`u $B_\rho(0,r)=\{x\in
X_\rho \, /\, \rho(x) < r\}$.

\begin{prop} 1) la famille ${\cal B}$ est une base de filtres.

\quad \quad \quad \quad \quad \quad 2) Tout \'el\'ement de ${\cal B}$ est \'equilibr\'e, absorbant et si de plus $\rho$ est un convexe alors tout \'el\'ement de ${\cal B}$ est convexe.
\end{prop}
{\bf Preuve de la proposition2.1.}

1) ${\cal B}$ est une base de filtres. En effet:

 a) $\emptyset \notin {\cal B}$ car tout $B_\rho(0,r)$ contient son
centre.

b) Soient $B_\rho(0,r_1)$ et $B_\rho(0,r_2)$ dans ${\cal B}$, posons $r= \inf (r_1,r_2)$, on obtient

$$B_\rho(0,r) \subset B_\rho(0,r_1)\cap B_\rho(0,r_2).$$
 En effet, pour
 $z\in B_\rho(0,r)$ on a:
    \begin{eqnarray*}
\left \{
\begin{array}{c}
\rho(z) < r\leq r_1\\ \rho(z)< r\leq r_2
\end{array}
\right.
\end{eqnarray*}
 donc $z\in B_\rho(0,r_1)\cap B_\rho(0,r_2)$.\\
Par suite il existe $B_\rho(0,r)\in {\cal B}$ tel que
$B_\rho(0,r)\subset B_\rho(0,r_1)\cap B_\rho(0,r_2)$. D'o\`u
${\cal B}$ est une base de filtres.

 2) Soit $B_\rho(0,r)\in {\cal B}$.

i) $B_\rho(0,r)$ est \'equilibr\'e. En effet: \\
 Soit $\alpha\in \C$ tel que
$\lambda=|\alpha|\leq 1$, alors il existe $\theta \in \R$ tel que
$\alpha=\lambda e^{i\theta}$, prenons $x\in B_\rho(0,r)$ \\
on a:
$$\rho(\alpha x)=\rho(\lambda e^{i\theta}x)=\rho(\lambda x)\leq
\rho(x)< r ,$$
 ceci  impliquera que $\alpha x\in B_\rho(0,r)$.

ii) $B_\rho(0,r)$ est absorbant. En effet:

 Soit $x\in X_\rho$ donc
$\lim\limits_{\lambda \to 0}\rho(\lambda x)=0$ ce qui est
\'equivalent \`a
\\ $\forall r>0, \exists \delta >0,$ tel que $ 0< \lambda
<\delta$ on a $\rho(\lambda x)< r$, donc il existe $\lambda
>0$  tel que $\lambda x\in B_\rho(0,r)$ ce qui montre que
$B_\rho(0,r)$ est absorbant.\\

 Si maintenant $\rho$ est convexe,
soient $B_\rho(0,r)\in {\cal B}$,  $x,y \in
B_\rho(0,r)$ et $\lambda \in [0,1]$ on a:
 $$\rho(\lambda x + (1-\lambda)y)
\leq \lambda \rho(x)+(1-\lambda)\rho(y)< r ,$$
donc
$$\lambda x + (1-\lambda)y\in B_\rho(0,r).$$
 D'o\`u $B_\rho(0,r)$ est convexe.
\begin{theo}

 Supposons que le modulaire
$\rho$ v\'erifie la condition $\Delta_2$. Alors
 $X_\rho$ est un espace vectoriel topologique s\'epar\'e.
\end{theo}
{\bf Remarque.}

Rappelons que si $\rho$ v\'erifie la condition $\Delta_2 $, J.Musielak a montr\'e dans \cite{jm} que $X_\rho $ est un espace vectoriel topologique sans mention\'e la s\'eparation de la topologie. \\
 Nous cherchons gr\^ace \`a
la forme (3)  de la condition $\Delta_2$ une autre
d\'emarche (la topologie d\'efinit par une base de filtres ) de montrer ce th\'eor\`eme en pr\'ecisant que la topologie
est s\'epar\'ee.  D'o\`u l'int\'er\^et de la forme (3)  de la
condition $\Delta_2$.\\
{\bf  Preuve du Th\'eor\`eme 2.1.}

 La d\'emonstration est bas\'ee sur le Lemme 1.1 et
la Proposition 2.1.\\
 D'apr\`es la Proposition 2.1, la
famille ${\cal B}$ est une base de filtres, de plus tout \'el\'ement de
${\cal B}$ est \'equilibr\'e et absorbant. \\
D'autre part, pour tout
$B_\rho(0,r)$, il existe $\delta_0>0$ tel que
$$B_\rho(0,\delta_0) + B_\rho(0,\delta_0)\subset B_\rho(0,r).$$
 En effet, soit $r> \epsilon >0$, comme $\rho$ v\'erifie la condition $\Delta_2$, donc  $\exists
L>0, \exists \delta >0$, tel que si $\rho(x)<L$,
$\rho(x-y)<\delta$ alors $|\rho(y)-\rho(x)|<\epsilon$ (voir le Lemme1.1).\\
Posons
$$\delta_0=\inf(r-\epsilon ,L, \delta)$$
 et soit \\ $z\in
B_\rho(0,\delta_0)+B_\rho(0,\delta_0)$ donc $z=x+y$ avec
\begin{eqnarray*}
\left \{ \begin{array}{c} \rho(x)< \delta_0\\ \rho(y)<
\delta_0.
\end{array}
\right.
\end{eqnarray*}
On a: $y=z-x\in B_\rho(0,\delta_0)$ implique $ \rho(z-x)<
\delta_0\leq \delta $ et $\rho(x)< \delta_0\leq L$, \\
d'o\`u
$$\rho(z)< \epsilon + \rho(x) < \epsilon +\delta_0\leq
\epsilon +r -\epsilon =r.$$
 Ce qui montre que $z\in B_\rho(0,r)$.
 Par suite $B_\rho(0,\delta_0)+ B_\rho(0,\delta_0)\subset B_\rho(0,r)$.\\
  Par cons\'equent $X_\rho$ est un espace vectoriel topologique.
 La topologie est d\'efinie
par:  $${\cal T} =\{ G\ne \emptyset, G\subset X_\rho \, /\, \mbox{\,
si }\,  x\in G, \, \mbox{ \, alors }\, \exists V\in {\cal B} \,
\mbox{\, tel que }\, x +V\subset G\}\cup \{\emptyset \} .$$
 Il est facile de v\'erifier que ${\cal T}$ est une
topologie de $X_\rho$. Il reste \`a montrer que $(X_\rho,{\cal T}
)$ est s\'epar\'e. Pour cel\`a, soient $x,y$ dans $X_\rho$ tel que $x\ne y$. \\
Supposons que pour tout $V_x$ voisinage de $x$ et $V_y$ voisinage
de $y$ on a: $V_x\cap V_y\ne \emptyset $.\\ Prenons $z\in (x +
B_\rho(0,\frac 1n))\cap (y + B_\rho(0,\frac 1n))$, $n\in \N^{*}$.
On a:
\begin{eqnarray*}
\left \{ \begin{array}{c} \rho(x-z)< \frac 1n\\ \rho(y-z)<
\frac 1n.
\end{array}
\right.
\end{eqnarray*}
  Comme $\rho$ v\'erifie la condition $\Delta_2$, alors d'apr\`es le Lemme1.1,
 pour tout  $ \epsilon >0$, ils existent $ L>0$ et $\delta >0$, tels que si
$\rho(x)<L$ et $\rho(y-x)<\delta$ alors $|\rho(y)-\rho(x)|< \frac
\epsilon 2$. \\
 On pose  $Y = y-x$  et $X = z-x$.  Alors on a:
\begin{eqnarray*}
\left \{ \begin{array}{c} \rho(X)=\rho(x-z)< \frac 1n\\
\rho(Y-X)=\rho(y-z)< \frac 1n .
\end{array}
\right.
\end{eqnarray*}
Pour $n$ assez grand tel que $\frac 1n \leq \inf(L,\delta,
\frac \epsilon 2)$, et par le Lemme1.1 on a:
 $$ \rho(Y) = \rho(y-x)< \rho(z-x) +\frac
\epsilon 2\leq \frac \epsilon 2+\frac \epsilon 2=\epsilon .$$
 Donc pour tout $\epsilon >0$, on a: $\rho(y-x)< \epsilon
$, ce qui montre que $ \rho(x-y) =0$, d'o\`u  $x = y$ ce qui contredit le
fait que $x\ne y$, par suite il existe un voisinage $V_x$ de $x$
et un voisinage $V_y$ de $y$ tel que $V_x\cap V_y=\emptyset $.

{\bf Remarque.}\\
Si de plus le modulaire $\rho$ est convexe alors $X_\rho $ est un espace vectoriel topologique localemment convexe. En effet:\\
 D'apr\`es le Th\'eor\`eme 2.1, $X_\rho$ est un espace
vectoriel topologique s\'epar\'e, et d'apr\`es
la Proposition 2.1 (2) tout \'el\'ement de ${\cal B}$ est convexe,
d'o\`u $X_\rho$  est un espace vectoriel topologique localement
convexe s\'epar\'e.
\section{ Caract\'erisation des $\tau $-ferm\'es de l'espace $X_\rho$.}

Rappelons que dans la th\'eorie des espaces modulaires, on trouve  les   s\'equentiellement ferm\'es voir par exemple page 19 de \cite{jm}. Dans ce travail, nous pr\'esentons la d\'efinition et la caract\'erisation des $\tau$-ferm\'es de l'espace $X_\rho .$
\begin{deff} Soient $\rho$ un modulaire v\'erifiant la condition $\Delta_2$ et $F \subset X_\rho$, on dit que $F$ est un ferm\'e pour la
topologie $\tau$ ( $\tau$-ferm\'e) si et seulement si son compl\'ementaire dans $X_\rho$, qu'on notera $C^F$, est un
ouvert pour la topologie $\tau$.
\end{deff}
\begin{theo} Soient $\rho$ un modulaire v\'erifiant la condition $\Delta_2$ et  $F\subset X_\rho$, alors on a:
\\ $F$
est un ferm\'e au sens de la topologie $\tau$  si et seulement si
pour toute suite $(x_n)_{n\in \N}\subset F$ tel que $x_n
\stackrel{\rho}{\longrightarrow} x$, alors $x\in F$.
\end{theo}

Pour montrer le Th\'eor\`eme 3.1  nous avons besoin de la proposition suivante.
\begin{prop}
Soient $\rho$ un modulaire v\'erifiant la condition $\Delta_2$ et  $F$ un ferm\'e au sens de la topologie $\tau $ dans $X_\rho$, alors
$$  x\in F \Longleftrightarrow \forall \epsilon >0 , B_\rho (x,\epsilon ) \cap F \ne \emptyset $$
\end{prop}
{\bf Preuve de la proposition.}

Soit $x\in X_\rho $,
\begin{eqnarray*}
x\not\in F & \Longleftrightarrow & x\in C^F \, \mbox{\,  qui est un ouvert
pour la topologie } \tau \\
 & \Longleftrightarrow & \exists B_\rho(0,\epsilon )\in
{\cal B}\, / x+ B_\rho(0,\epsilon )=B_\rho(x,\epsilon )\subset
C^F\\
& \Longleftrightarrow &  \exists \epsilon >0, \, \mbox{\, tel que \, }
B_\rho(x,\epsilon )\cap F=\emptyset
\end{eqnarray*}
D'o\`u $x\in F \Longleftrightarrow  \forall \epsilon >0, B_\rho(x,\epsilon
)\cap F\ne\emptyset $.\\

{\bf Preuve du Th\'eor\`eme 3.1.}

Soient $F$ un ferm\'e pour la topologie $\tau$ et $(x_n)_{n\in \N}$ une suite dans  $ F$  tel que $x_n
\stackrel{\rho}{\longrightarrow} x$ alors  pour tout $\epsilon>0 $, il existe $n_0\in \N$ tel que si $n>n_0$, on a $x_n\in B(x,\epsilon )$,
ceci impliquera que
$$\forall \epsilon >0 \ ,   B(x,\epsilon )\cap F\ne \emptyset ,$$
 d'o\`u  d'apr\`es la Proposition 3.1, $x\in F $.\\

Inversement, supposons que $F$ n'est pas un ferm\'e pour la topologie $\tau$,
donc $C^F$ n'est pas un ouvert pour la topologie $\tau$, ce qui
entraine qu'il existe $x\in C^F$ tel que pour tout $\epsilon >0$,
$B_\rho(x,\epsilon)\not\subset C^F$ d'o\`u, il existe $x\in C^F$
tel que pour tout $\epsilon >0$, $B(x,\epsilon)\cap F\ne \emptyset $.
Donc,\\
 pour $\epsilon =1$ il
existe $x_1\in B_\rho(x,1)\cap F$\\
 pour $\epsilon =\frac 12$ il
existe $x_2\in B_\rho(x,\frac 12)\cap F$
$$.\\ . \\ . $$
\noindent pour $\epsilon =\frac 1n$ il existe $x_n\in
B_\rho(x,\frac 1n)\cap F$.\\
 On obtient donc la suite $(x_n)_{n\in \N}\subset F$ tel que
  $x_n \stackrel{\rho}{\longrightarrow} x$,
   alors $x\in F$, contradiction car $x\in C^F$. D'o\`u $F$
est un ferm\'e pour la topologie $\tau$.
\section{ Certains r\'esultats de s\'eparations dans les espaces modulaires.}
\quad Consid\`erons $\rho$ un modulaire   v\'erifiant la condition
$\Delta_2$ et soit $A$ un sous-ensemble de $X_\rho$, on pose par
d\'efinition pour tout $x\in X_\rho$,
 $$\rho(x,A)=\inf\{\rho(x-y),\, y\in A\}.$$

Nous pr\'esentons un r\'esultat de s\'eparation dans les espaces
modulaires.

\begin{theo} Soient $\rho$ un modulaire  v\'erifiant la condition
$\Delta_2$, $A $ une partie  ferm\'ee de $X_\rho$ et $x_0\notin A$.
Alors il existe un voisinage ouvert  $V_{x_0}$ de $x_0$ tel que
$V_{x_0}\cap A =\emptyset $.
\end{theo}
Pour montrer le Th\'eor\`eme 4.1, nous avons besoin du r\'esultat suivant.

\begin{prop} Soient $\rho$ un modulaire v\'erifiant la condition $\Delta_2$ et $A \subset X_\rho$ on a:\\
 $\rho(x,A)=0 $ si et seulement si $ x\in {{\overline{A}}^{\rho}} $, o\`u ${{\overline{A}}^{\rho}} $ est la fermuture de $A$ pour la topologie $\tau$.
\end{prop}
{\bf Preuve.}

On a:
$$\rho(x,A)=\inf\{\rho(x-y),\, y\in A\}=0 ,$$
 donc pour tout
$\epsilon =\frac 1n$,  il existe $y_n\in A$  tel que
$\rho(x-y_n)<\frac 1n$ ce qui montre qu'il existe une suite
$(y_n)_{n\in \N}\subset A$  tel que $y_n
\stackrel{\rho}{\longrightarrow} x$ d'o\`u $x\in
{\overline{A}^\rho}$.

Inversement, soit $x\in {{\overline{A}}^\rho}$, on a d'apr\`es  le Th\'eor\`eme 3.1, il existe une
suite $(y_n)_{n\in \N}\subset A$  tel que $y_n
\stackrel{\rho}{\longrightarrow} x$, donc pour tout $\epsilon >0$
il existe $n_0$, pour $n>n_0$ on a:
 $$\rho(x,A)\leq \rho(x-y_n)<\epsilon.$$
 D'o\`u
$$\rho(x,A)=0.$$
{\bf Preuve du Th\'eor\`eme 4.1.}

 On a d'apr\`es la Proposition 4.1, $x_0\notin A $ si seulement si $ \rho(x_0 ,A)=r>0$.\\
Comme $\rho$ v\'erifie la condition $\Delta_2$ alors d'apr\`es le Lemme1.1, pour
$\epsilon =\frac r3>0$, ils existent $L>0, \delta >0$ si $\rho(x)<L$ et
$\rho(y-x)<\delta$ alors $|\rho(y)-\rho(x)|<\epsilon $. \\
De plus il
existe un $m_0\in \N^{*}$ tel que $\frac rm<\inf (L,\delta)$ d\`es
que $m>m_0$, on choisit $ m_1 \geq max (3,m_0 )$ et on consid\`ere le voisinage ouvert de $x_0$:
$$V_{x_0}=x_0+B_\rho(0,\frac r{m_1}).$$
 Supposons qu'il existe $y\in V_{x_0}\cap A$, comme $A$ est ferm\'e d'apr\`es la Proposition 4.1, il existe
une suite $(y_n)_{n\in \N}\subset A$  tel que
$y_n\stackrel{\rho}{\longrightarrow} y$. \\
 Posons $X_n=y-y_n$ et
$Y_n=x_0-y_n$ comme $y_n\in A$ et $x_0\notin A$ alors
$\rho(Y_n)\geq r$. \\
D'autre part,
$$\rho(X_n)=\rho(y-y_n)< \frac r{m_1}<\inf (L,\delta),$$
 d\`es que $n>n_0$ et
$$\rho(X_n-Y_n)=\rho(x_0-y)< \frac r{m_1}<\inf (L,\delta),$$
 d'o\`u
$r\leq \rho(Y_n) < \rho(y-y_n)+\epsilon \leq \frac{r}{m_1} +\frac
{r}{3}\leq \frac{2r}{3}$ d\`es que $n>n_0$ absurde donc
$V_{x_0}\cap A=\emptyset $.\\

{\bf Remarque.}

 Si de plus $\rho$ v\'erifie la propri\'et\'e de Fatou on aura:\\
$\overline{B(0,r)}= \overline{B_\rho }(0,r)=\{x\in X_\rho\, /\, \rho(x)\leq
r\}$ est une boule ferm\'ee pour la topologie $\tau$. Toute boule ferm\'ee de centre $x$ et de rayon $r$ sera not\'ee $B_f (x,r) $ ( voir \cite{krk} ).
\begin{coro}
Sous les m\^emes hypoth\`eses du th\'eor\`eme pr\'ec\'edent et si
de plus $\rho$ v\'erifie la propri\'et\'e de Fatou alors
${\overline{V_{x_0}}}^\rho\cap A=\emptyset $.
\end{coro}
{\bf Preuve.}

D'apr\`es le Th\'eor\`eme 4.1, il existe
$V_{x_0}=x_0+B_\rho(0,\frac r{m_1})$ tel que $V_{x_0}\cap
A=\emptyset $. \\
 D'abord on montre que ${\overline{V_{x_0}}}^\rho =x_0+B_f(0,\frac r{m_1})$. En effet:\\
 Soit $y\in {\overline{V_{x_0}}}^\rho$ alors d'apr\`es la Proposition 4.1, il existe une suite
$(y_n)_{n\in\N}\subset B_f (0,\frac r{m_1})$  tel que
$$x_0+y_n \stackrel{\rho}{\longrightarrow} y,$$
 ceci  impliquera que  $y_n
\stackrel{\rho}{\longrightarrow} y-x_0$. En effet: soient
$$Y_n=y_n-(y-x_0) \ \ {\mbox et}\ \  X_n=2(y_n-(y-x_0)),$$
 il est facil de voir que
$Y_n \stackrel{\rho}{\longrightarrow} 0$  et comme $\rho $ v\'erifie la condition $\Delta_2 $ alors
$X_n \stackrel{\rho}{\longrightarrow} 0$.\\
 D'o\`u pour $\epsilon >0, \exists L>0, \exists
\delta>0$ tels que
$$\rho(X_n) < \inf (L,\delta ,\frac \epsilon 2),$$
 et
$$\rho(Y_n-X_n)= \rho(Y_n ) <\inf (L,\delta ,\frac \epsilon 2),$$
 d\`es que $n\geq n_0$, donc
$$\rho(Y_n)=\rho(y_n-(y-x_0)) < \inf (L,\delta ,\frac
\epsilon 2)+ \frac \epsilon 2\leq \frac \epsilon 2 +\frac \epsilon
2= \epsilon ,$$
 d\`es que $n\geq n_0$, par suite
$$y_n \stackrel{\rho}{\longrightarrow} y-x_0 \in
{\overline{B_\rho (0,\frac r{m_1})}}^\rho =B_f (0,\frac r{m_1}),$$
 donc
$$y=x_0+(y-x_0)\in x_0 +B_f(0,\frac r{m_1}),$$
 d'o\`u
$${\overline{V_{x_0}}}^\rho \subset x_0+B_f(0,\frac r{m_1}).$$
Inversement, soit
$$x_0+y\in  x_0 +B_f(0,\frac r{m_1}),$$
 d'apr\`es la Proposition 4.1, il
existe $(y_n)_{n\in\N}\subset B_\rho (0,\frac r{m_1})$  tel que $y_n
\stackrel{\rho}{\longrightarrow} y$ et on
a aussi la suite $(x_0+y_n)_{n\in\N}\subset V_{x_0}$ v\'erifiant
$x_0+y_n \stackrel{\rho}{\longrightarrow} x_0+y$ d'o\`u
$$x_0+y\in {\overline{V_{x_0}}}^\rho . $$
 Pour terminer, on
suit les m\^emes d\'emarches de la preuve du Th\'eor\`eme 4.1,
on  aura:
$${\overline{V_{x_0}}}^\rho \cap A=\emptyset . $$
\begin{theo} Soit $\rho$ un  modulaire convexe  v\'erifiant la condition
$\Delta_2$. Alors l'espace modulaire
$X_\rho$ est  normal.
\end{theo}
{\bf Preuve du Th\'eor\`eme 4.2.}\\

Le modulaire $\rho$ est convexe et v\'erifie la condition $\Delta_2$, alors les convergences au sens du modulaire et  au sens de la norme sont \'equivalente, donc  l'espace modulaire $X_\rho$ est maîtrisable, alors $X_\rho$ est normal.

{\bf Adresse:} D\'{e}partement de Math\'{e}matiques et
Informatique.\\ Fac. des Sc. Univ.
Mohammed V-Agdal \\ B.P 1014, Rabat, Maroc. \\

E-mail: hajid32@yahoo.fr  \\ \\ \\

{\bf Remerciement.}\\

Je tiens \`a remercier le Professeur E.Hanebaly  pour ses remarques sur ce travail.

\end{document}